\documentclass[]{article}
\usepackage{amsfonts}
\usepackage{dsfont}
\usepackage{amsmath}
\usepackage{verbatim}
\usepackage{amsthm}
\usepackage{enumerate}
\usepackage[affil-it]{authblk}

\theoremstyle{plain}
\newtheorem{definition}{Definition}
\newtheorem{theorem}{Theorem}

\newtheorem{lemma}{Lemma}

\theoremstyle{definition}
\newtheorem{remark}{Remark}

\title{Limits of random trees\thanks{MSC2010 Subject Classification: 05C80}}
\author{Attila De\'ak\\
\texttt{deak51@cs.elte.hu}}
\affil{MTA-ELTE "Numerical Analysis and Large Networks" Research Group, 1117 Budapest, P\'azm\'any P. s. 1/c., Hungary}

\begin{document}
\maketitle
\let\thefootnote\relax\footnote{Research was supported by the ERC Grant Nr.: 227701}
\let\thefootnote\relax\footnote{The final publication is available at http://link.springer.com/article/10.1007\%2Fs10474-013-0321-0}

\begin{abstract}
Local convergence of bounded degree graphs was introduced by Benjamini and Schramm \cite{B-Sch_RoDLoFPG}.
This result was extended further by Lyons \cite{Lyons_AEoST} to bounded average degree graphs.
In this paper, we study the convergence of a random tree sequence $(T_n)$, where the probability of a given tree $T$ is proportional to $\prod_{v_i\in V(T)}d(v_i)!$.
We show that this sequence is convergent and describe the limit object, which is a random infinite rooted tree.\\

\end{abstract}

\bigskip


\section{Introduction}


Limits of graph sequences with bounded degree have been studied extensively over the last decade.
A natural extension is to study the case when we only require bounded average degree.
For instance trees have average degree less than $2$, but in general a sequence of trees can have unbounded maximum degree.
A limit theory for trees has been established by Aldous \cite{Aldous_CRT} and also by Elek and Tardos \cite{E-T_LoT}.

We do not follow their path, but use the limit theory for bounded average degree graphs, described by Lyons \cite{Lyons_AEoST}.
Define $T_n$ as the random tree on the nodes $\{1,2,\cdots , n\}$ so that for a given tree $T$ with degrees $d_i$ we have
$$
\mathds{P}(T_n=T)={\prod_{i=1}^nd_i!\over (n-2)!{3n-3 \choose n-2}}.
$$
We will show that $T_n$ converges and has a limit, a random infinite rooted tree. Let
${\cal A}_n$ denote the set of trees with $n$ nodes. For motivation consider the following process on ${\cal A}_n$:

\begin{itemize}
\item
Choose a random edge and also one of its endpoints uniformly $(X, V_{old})$.
\item
Take a uniformly chosen neighbor of $V_{old}$: $V_{new}$.
\item
If $X=V_{new}$, then do nothing, or else remove the edge $(X, V_{old})$ and add a new edge $(X, V_{new})$.
\end{itemize}

This clearly defines a Markov chain $(A_t^n)$ on ${\cal A}_n$. Let $\Pi_n$ denote the stationary distribution of the process defined above. 
It is easy to prove that the Markov chain defined above is reversible and aperiodic. From the reversibility we can also compute the stationary distribution:
$$\Pi_n(T)=\frac{\prod_{i=1}^n d_i!}{C},$$
where $C=(n-2)!{3n-3\choose n-2}$, see Remark \ref{konst} below.

The distribution of $T_n$ depends only on the degrees of its vertices and not on the exact structure of the tree itself. The natural question is: Considering only trees, can the degree sequence alone determine the convergence of a tree sequence? We show that with the above distribution the random tree sequence converges. We do not know how much this result can be extended by defining other distributions, which ensure random local convergence.

Our paper is organized as follows.
In Section \ref{basic_def} we give the basic definitions and state the main theorem of this paper.
Then in Section \ref{deg_dist_max_deg} we prove basic bounds and asymptotics for the degrees of $T_n$, such as the expected number of degree $d$ vertices and the expected value of the maximum degree.
In Section \ref{lab_subg_dens} we investigate subgraph densities and in Section \ref{limit} we describe the limit object.

\section{Basic definitions}
\label{basic_def}

For a finite simple graph $G$, let $B_G(v,R)$ be the rooted $R$-ball around the node $v$, that is the subgraph induced by the nodes at distance at most $R$ from $v$.
Given a positive integer $R$, a finite rooted graph and a probability distribution $\rho$ on rooted graphs, let $p(R,F,\rho)$ denote the probability that the graph $F$ is rooted isomorphic to the $R$-ball around the root of a rooted graph chosen with distribution $\rho$.
For a finite graph $G$, let $U(G)$ denote the distribution of rooted graphs obtained by choosing a uniform random node of $G$ as the root.

\begin{definition}
Let $(G_n)$ be a sequence of random finite graphs, $\rho$ a probability distribution on rooted graphs.
We say that the random local limit of $G_n$ is $\rho$, if for any positive integer $R$ and finite rooted graph $F$, we have
$$
\lim_{n\rightarrow \infty}\mathds{P}(|p(R,F,U(G_n))-p(R,F,\rho)|>\epsilon)=0.
$$
\end{definition}

\begin{theorem}
\label{tree_limit}
Let $X_n$ be a random tree from the distribution $\Pi_n$. $X_n$ has a random local limit, which is an infinite rooted random tree.
\end{theorem}

Let $T$ be a tree on $n$ nodes and denote by $X_n$ a random tree with distribution $\Pi_n$. We know that
$$
\mathds{P}(X_n=T)={1\over C}\prod_{i=1}^nd_i!,
$$
where $d_i$ is the degree of the $i$th vertex of the tree $T$. It is easy to see that there are ${n-2 \choose d_1-1,d_2-1,\cdots ,d_n-1}$ trees that realize the same degree sequence.
From this it follows that for a given degree sequence $(d_i)_{i=1}^n$, we get
\begin{multline}
\label{deg_seq_prob}
\mathds{P}((d_i)_{i=1}^n )=\frac{1}{C}\prod_{i=1}^n d_i!{n-2 \choose d_1-1,d_2-1,\cdots ,d_n-1}=\\
=\frac{1}{C}\prod_{i=1}^n d_i!\frac{(n-2)!}{\prod_{i=1}^n (d_i-1)!}=\frac{(n-2)!}{C}\prod_{i=1}^n d_i.
\end{multline}
To be able to compute probabilities about the degree sequence, we need to calculate the sum of the product of possible degree sequences.
The following lemma states an easy result about this.
\begin{lemma}
\label{sum_prod_deg}
$$\sum_{\substack{\sum_{i=1}^n d_i=m\\d_i>0}}\prod_{i=1}^n d_i={n+m-1 \choose m-n}.$$
\label{alt_prod_osszeg}
\end{lemma}

\noindent
{\bf Proof.}
Let $f(x)=\left (\sum_{i=1}^\infty ix^i\right)^n$ and denote the sum in the lemma by $M$. It is easy to see that $M$ is the coefficient
of $x^m$ in:
$$
f(x)=\left (\sum_{i=1}^\infty ix^i\right)^n=\left ({x\over (1-x)^2}\right)^n={x^n\over (1-x)^{2n}},
$$
hence it follows that $M=(-1)^{m-n}{-2n \choose m-n}={m+n-1 \choose m-n}$.
\qed

\begin{remark}
\label{konst}

From Lemma \ref{alt_prod_osszeg}, it follows that $\sum_{(d_i)} \prod_i d_i={3n-3 \choose n-2}$,
where the summation is over the possible degree sequences.
From this using (\ref{deg_seq_prob}) $C=(n-2)!{3n-3 \choose n-2}$ also follows.
\end{remark}

\section{Degree distribution and the maximum degree}
\label{deg_dist_max_deg}
%
Using Lemma \ref{sum_prod_deg} we can compute the expected number of degree $d$ vertices and the maximum degree.
\begin{lemma}
\label{lemma_fokszam_eloszlas}
Let $M=\sum_{i=1}^k x_i$ and assume that $M \leq n+k-2$.
Then the probability that the $i$th vertex has degree $x_i$ ($i=1,\ldots ,k$) is
$$\mathds{P}(d_i=x_i; \,i=1,\cdots,k)=\frac{{3n-M-k-3\choose n-M+k-2}}{{3n-3\choose n-2}}\prod_{i=1}^k x_i.
$$
\end{lemma}

\noindent
{\bf Proof.}
Using the above notions with (\ref{deg_seq_prob}) we get:
$${3n-3\choose n-2}\mathds{P}(d_i=x_i \,;i=1,\cdots ,k)=\sum_{\substack{\sum d_i=2n-2 \\ d_j=x_j\, j\leq k}}\prod_i d_i=$$
$$\prod_{i=1}^k x_i\sum_{\substack{\sum_{j=k+1}^n d_j= \\ =2n-M-2}}\prod_{j=k+1}^n d_j={3n-M-k-3\choose n-M+k-2}\prod_{i=1}^k x_i,$$
where the last equation follows from Lemma \ref{alt_prod_osszeg} and this is what we wanted to prove.
\qed

Let us denote the maximum degree in a tree $T$ by $D(T)=\max_i d_i$.
Now the following theorem is true:

\begin{theorem}
For every $\epsilon>0$, we have
$$(1-\epsilon-o(1))\log_3n\leq \mathds{E}(D(X_n))
\leq (1+\epsilon+o(1))\log_3n, \textrm{ as } n\rightarrow \infty.$$
\end{theorem}
We bound the expected value by bounding the probabilities $\mathds{P}(D>k)$ and $\mathds{P}(D\leq k)$ with Lemma \ref{lemma_maximum_also_becsles} and Lemma \ref{lemma_maximum_felso_becsles} respectively.
To simplify the notation, we let $D=D(X_n)$.

\begin{lemma}\label{lemma_maximum_also_becsles}
For every $\delta>0$, there exists $n_0$, such that $\forall n>n_0$
$$\mathds{P}(D>k)<(1+\delta)\frac{4}{3}\frac{k}{3^{k}}n.$$
\end{lemma}

\noindent
{\bf Proof:}
From Lemma \ref{lemma_fokszam_eloszlas} it follows that
$$\mathds{P}(d_i=k)=k\frac{{3n-k-4 \choose n-k-1}}{{3n-3 \choose n-2}}< {4\over3}{(1+\delta)k\over3^k},\ \forall \delta>0,n>n_0(\delta),$$
and so
$$\mathds{P}(d_i> k)<\sum_{j=k+1}^n {4\over3}{(1+\delta)j\over3^j}<\sum_{j=k+1}^\infty {4\over3}{(1+\delta)j\over3^j}=
(1+\delta){2k+3\over 3^{k+1}}.$$
Now the statement of the lemma follows, as $\mathds{P}(D> k)<n\mathds{P}(d_i>k)$.
\qed

\begin{lemma}\label{lemma_maximum_felso_becsles}
$$\mathds{P}(D\leq k)<{1\over 9}\sqrt{n}\,\exp\left\{-n{4\over 9}{(k+1)\over 3^k}\right\}$$
\end{lemma}

\noindent
{\bf Proof:}
It is easy to see that
$$\mathds{P}(D\leq k)={1\over {3n-3 \choose n-2}}\sum_{\substack{d_i\leq k\\ \sum_i d_i=2n-2}}\prod_{i=1}^n d_i.$$
Here the sum is just the  coefficient of $x^{n-2}$ in  $P(x)=\Bigl (1+2x+3x^2+\ldots +kx^{k-1}\Bigr )^n$.
As $P(x)=\sum a_ix^i$ and $\forall i,\ a_i\geq 0$, we get an upper bound on $a_i$:
$$a_i\leq {P(x_0) \over x_0^i},$$
for any $x_0>0$.
Hence it follows that
$${3n-3 \choose n-2}\mathds{P}(D\leq k)\leq P({1\over3})3^{n-2}=3^{n-2}\left (\sum_{i=1}^k {i\over 3^{i-1}} \right)^n
\leq 3^{n-2}\left ({9\over 4}- {k+1 \over 3^{k}} \right)^n=$$
$$= 3^{n-2}\left({9\over 4}\right)^n \left(1-{4\over 9}{k+1 \over 3^k}\right)^n<
{1\over 9}\left({27\over 4}\right)^n \exp\left\{-n{4\over 9}{k+1\over 3^k}\right\}.$$
By using Stirling's formula, we get
$$\mathds{P}(D\leq k)<{1\over 9}\sqrt{n}\,\exp\left\{-n{4\over 9}{k+1\over 3^k}\right\},$$
the desired inequality.
\qed

\noindent
Now we can upper bound the expected value as follows:
$$\mathds{E}(D)\leq k_1\mathds{P}(D\leq k_1)+k_2\mathds{P}(D>k_1)+(n-1)\mathds{P}(D>k_2),$$
where
$$k_1= (1+\epsilon)\log_3 n \textrm{ and } k_2=\log^3_3 n,$$
Therefore it follows from Lemmas \ref{lemma_maximum_also_becsles} and \ref{lemma_maximum_felso_becsles} that
$$\mathds{E}(D)\leq (1+\epsilon)\log_3 n+{\log^3_3 n \over n^{\epsilon +o(1)}}+{\log^3_3 n\over n}=(1+\epsilon +o(1))\log_3 n,\textrm{ as } n\rightarrow\infty.$$
Also
$$\mathds{E}(D)\geq (1-\epsilon)\log_3 n \mathds{P}\left(D>(1-\epsilon)\log_3 n\right)\geq (1-\epsilon -o(1))\log_3 n,\textrm{ as } n\rightarrow \infty,$$
which  proves our theorem.
\qed


\section{Labeled subgraph densities}
\label{lab_subg_dens}

Form now on let $T$ be a fixed tree on $k$ nodes. Assign values $r_i$
to each node $i$ of $T$ and call it the remainder degree of node $i$. Let
$S=\{s_1,s_2,\ldots ,s_k\}$ be an ordered subset of $[n]=\{1,2,\ldots ,n\}$
$(s_i\neq s_j,\,i\neq j)$. Now
by $T_S$ we denote the tree $T$, with label $s_i$ at node $i$.

Similarly let
$F$ be a forest with $m_F$ nodes and $c_F$ connected components. Denote these
components by $C_1,\cdots C_{c_F}$. As above, define the remainder degrees for $F$, and denote them by $(r_1,\ldots r_{m_F})$.
Denote by $F_S$ the labeled forest with label $s_i$ at node $i$.

For two labeled trees $T_S,T'_{S'}$, with remainder degrees $r,r'$
we define the operation gluing in the usual way.
We identify nodes $i\in V(T),j\in V(T')$ if $s_i=s'_j$ and keep the label
$s_i(=s'_j)$. For nodes $i\ (j)$, where $s_i\notin S'(s'_j\notin S)$, we do
nothing. If $S\cap S'=
\emptyset$, then the resulting graph is just the disjoint union.
We say that a gluing is valid, if it results in a forest and
$s_i=s'_j\Rightarrow r_i=r'_j\,\forall i,j$, and denote it by
$g(T_S,T'_{S'})$. We define the gluing of two labeled forests similarly.
Let
$$X_S^T=I\{\phi :i\mapsto s_i \text{ is a homomorphism from } T \text{ to }
X_n\text{, and }\forall i,\ d_{s_i}=d^T_i+r_i\}$$
$$X_n^T=\sum_{S\subseteq [n]}X_S^T=inj(T,X_n).$$
We define $X^F_S,X^F_n$ similarly for a forest $F$.
Also let $A_S^T=X^T_S-\mathds{E}(X^T_S)$.
If it does not make any confusion, we use $m_1,c_1,X_S,A_S$
instead of $m_{F_1},c_{F_1},X_S^T,A_S^T$.

\begin{theorem}
\label{labeled_tree_deviation}
For an arbitrary tree $T$ on $k$ nodes with remainder degrees $(r_i)_{i=1}^k$ and $\forall \epsilon>0$, we have
$$\mathds{P}\left({|X^T_n-\mathds{E}(X^T_n)|\over n}>\epsilon\right)\leq
{c_{k,r} \over n^2\epsilon^4}.$$
\end{theorem}
We want to bound the deviation from the expectation by bounding the $4$th moment of $X_n^T$.

\begin{lemma}
\label{4th_moment_bound}
With the above notions
$$\mathds{E}\left(\left(X_n^T-\mathds{E}(X^T_n)\right)^4\right)=O(n^2).$$
\end{lemma}

\begin{remark}
Theorem \ref{labeled_tree_deviation} is a direct corollary of Lemma \ref{4th_moment_bound} as
\begin{equation*}
\mathds{P}\left({|X^T_n-\mathds{E}(X^T_n)|\over n}>\epsilon\right)\leq
{\mathds{E}\left(\left(X_n^T-\mathds{E}(X^T_n)\right)^4\right) \over
n^4\epsilon^4}=
{c_{k,r} \over n^2\epsilon^4}.
\end{equation*}
\end{remark}
To get this bound we regard $X_n^T$ as a sum of indicator variables $X_S^T$, so we need to compute only probabilities $\mathds{P}(X_S^T=1)$ and $\mathds{P}(X_S^F=1|X_{S'}^{F'}=1)$ for special forests $F,F'$.

\begin{lemma}
\label{forest_prob}
Let $F$ be an arbitrary forest on $m$ nodes with remainder degrees $r=(r_1,\cdots,r_m)$.
Let $R=\sum_ir_i$ and denote by $c$ the number of connected components of $F$.
The probability that on an ordered subset $S=(s_1,\cdots,s_m)$ of the nodes of $X_n$ we see the forest $F$ and $\forall i,\ d_{s_i}=d^F_i+r_i$ is
$$
\mathds{P}(X_S^F=1)=\frac{(n-m+c-2)!}{(n-2)!}{{3n-R-3m+2c-3\choose n-R-m+2c-2}\over {3n-3\choose n-2}}H(r,F),
$$
where $H(r,F)$ is a constant depending only on $F$ and $r$.
\end{lemma}

\noindent
{\bf Proof:}
First we want to compute $\mathds{P}(X_S^F=1\,\big|\,(d_i)_{i=1}^n)$. Given the degree sequence $(d_i)$, the distribution of $X_n$ is uniform on the possible trees realizing $(d_i)$.
It is easy to check (simply by contracting the connected components to one node, counting the trees and blowing back the components) that the number of trees, with degree sequence $(d_i)$ and having $F$ on the first $m$ vertices and having $r_i$ edges going out from the forest at the $i$th node is

\begin{equation}
\displaystyle\prod\limits_{i=1}^{c}\left [{R_i! \over \prod\limits_{j\in C_i}r_j!} \right ]{n-m+c_F-2 \choose R_1-1,\cdots, R_{c}-1,d_{m+1}-1, \cdots,d_n-1}.
\end{equation}
Here $R_i=\sum_{j\in C_i}r_j$.
Since the number of trees realizing $(d_i)$ is
${n-2 \choose d_1-1,d_2-1,\cdots,d_n-1}$, it follows that

\begin{multline*}
\mathds{P}(X_{\{1,\cdots,m\}}^F=1|d_i)=\\
={\displaystyle{n-m+c-2 \choose R_1-1,\cdots, R_{c}-1,d_{m+1}-1, \cdots,d_n-1}
\over \displaystyle {n-2 \choose d_1-1,d_2-1,\cdots,d_n-1}}\displaystyle\prod\limits_{i=1}^{c}
{R_i! \over \prod\limits_{j\in C_i}r_j!}.
\end{multline*}

After simplifying, and using symmetry, we get the following equation:
\begin{multline}
\label{prob_x_on_T_2} 
\mathds{P}(X_S^F=1|d_i)=\\
={(n-m+c-2)! \over (n-2)!}
\prod_{i=1}^{c}\left[ R_i \prod_{j\in C_i}(d_{j}^F+r_j-1)\cdot\ldots\cdot (r_j+1)\right].
\end{multline}

Now to get the desired probability, we use Lemma \ref{lemma_fokszam_eloszlas} with $k=m$ and $M=(2m-c)+R$:

\begin{multline*}
\mathds{P}(X_S^F=1)=\sum_{d_i}
\mathds{P}(X_S^F=1|d_i){\prod_{i=1}^n d_i \over {3n-3\choose n-2}}=\\
={(n-m+c-2)!\over (n-2)!}
{{3n-R-3m+2c-3\choose n-R-m+2c-2}\over {3n-3\choose n-2}}
\prod_{i=1}^{c}\left[R_i\prod_{j\in C_i}\left[(d_j^F+r_j)\cdots (r_j+1)\right]\right]=\\
={(n-m+c-2)!\over (n-2)!}{{3n-R-3m+2c-3\choose n-R-m+2c-2}\over {3n-3\choose n-2}}H(r,F),
\end{multline*}
which is the claimed equation.
\qed

\begin{lemma}
\label{forest_cond_prob}
Let $F_1,F_2$ be two forests with remainder degrees $r,r'$ and labels $S_1,S_2$.
If there is a valid gluing of the labeled forests $F_1$ and $F_2$ then let $F_{1,2}=g(F_1,F_2)$.
Also let $m_{1,2}=|V(F_{1,2})|$, $c_{1,2}=\{$the number of components of $F_{1,2}\}$, $R_{1,2}=\sum_{V(F_{1,2})}r_i$. We have
\begin{multline}
\label{cond_forest_prob}
\mathds{P}(X_{S_1}^{F_1}=1|X_{S_2}^{F_2}=1)=\\
{(n-m_{1,2}+c_{1,2}-2)! \over (n-m_2+c_2-2)!}
{{3n-R_{1,2}-3m_{1,2}+2c_{1,2}-3\choose n-R_{1,2}-m_{1,2}+2c_{1,2}-2}
\over
{3n-R_2-3m_2+2c_2-3\choose n-R_2-m_2+2c_2-2}}
{H(r_{1,2},F_{1,2})\over H(r',F_2)},
\end{multline}
\end{lemma}

\noindent
{\bf Proof:}
The proof follows in a straightforward way from the definition of the conditional probability.
\qed

\begin{remark}
It is easy to see that ${H(r_{1,2},F_{1,2})\over H(r',F_2)}=H(r,F_1)$ if $S_1\cap S_2=\emptyset$.
\end{remark}

\begin{remark}
We can rewrite equation (\ref{cond_forest_prob}) as a rational polynomial in $n$ as follows.
\begin{multline}
\label{expected_value_as_polynom}
\mathds{P}(X_{S_1}^{F_1}=1|X_{S_2}^{F_2}=1){H(r',F_2) \over H(r_{1,2},F_{1,2})}=\\
={(n-m_{1,2}+c_{1,2}-2)!\over (n-m_2+c_2-2)!}{{3n-R_{1,2}-3m_{1,2}+2c_{1,2}-3\choose n-R_{1,2}-m_{1,2}+c_{1,2}-2}\over
{3n-R_2-3m_2+2c_2-3\choose n-R_2-m_2+c_2-2}}=\\
{2^{2(m_{1,2}-m_2)}\over 3^{R_{1,2}-R_2+3(m_{1,2}-m_2)-
2(c_{1,2}-c_2)}}
{P(n;R_{1,2},m_{1,2},c_{1,2},R_2,m_2,c_2)\over
Q(n;R_{1,2},m_{1,2},c_{1,2},R_2,m_2,c_2)}
\end{multline}
Where
\begin{multline*}
P(n;x_1,x_2,x_3,x_4,x_5,x_6)=\\
=n^{x_1-x_4+3(x_2-x_5)-2(x_3-x_6)}
-n^{x_1-x_4+3(x_2-x_5)-2(x_3-x_6)-1}\\
\left[
(x_1+x_4+x_2+x_5-x_3-x_6+3)
{(x_1-x_4+x_2-x_5-x_3+x_6)\over 2}+\right.\\
\left.+(2(x_2+x_5)-x_3-x_6+1){2(x_2-x_5)-x_3+x_6\over 4} \right ]
+O(n^{x_1-x_4+3(x_2-x_5)-2(x_3-x_6)-2})
\end{multline*}
and
\begin{multline*}
Q(n;x_1,x_2,x_3,x_4,x_5,x_6)=\\
=n^{x_1-x_4+4(x_2-x_5)-3(x_3-x_6)}-n^{x_1-x_4+4(x_2-x_5)-3(x_3-x_6)-1}\\
\left[(x_2+x_5-x_3-x_6+3){(x_2-x_5-x_3+x_6)\over 2}+\right.\\
\left.+(x_1+x_4+3(x_2+x_5)-2(x_3+x_6)+5){x_1-x_4+3(x_2-x_5)-2(x_3-x_6)
\over 6}\right]+\\
+O(n^{x_1-x_4+4(x_2-x_5)-3(x_3-x_6)-2}).
\end{multline*}
From this it also follows that
\begin{equation}
\label{cond_exp_order}
\mathds{E}(X_{S_1}^{F_1}|X_{S_2}^{F_2}=1)=
O(n^{m_2-m_{1,2}+c_{1,2}-c_2})
\end{equation}
and
\begin{equation}
\label{cond_exp_diff_order}
\mathds{E}(X_{S_1}^{F_1}|X_{S_2}^{F_2}=1)-\mathds{E}(X_{S_1}^{F_1})
=O(n^{m_2-m_{1,2}+c_{1,2}-c_2-1})=O(n^{-m_1+c_1-1}),
\end{equation}
if $S_1\cap S_2=\emptyset$ (since in this case $m_{1,2}-m_2=m_1$ and $c_{1,2}-c_2=c_1$).

\end{remark}

Now we have arrived at the main lemma of Chapter 4. Here we want to bound the following sum
$$
\sum_{S_1,S_2,S_3,S_4}\mathds{E}(A^T_{S_1}A^T_{S_2}A^T_{S_3}A^T_{S_4}).
$$

We split this sum in 5 parts, depending on the sizes of the intersections $S_i\cap S_j,\ (1\leq i<j \leq 4)$. We can change the indeces so that one of the following holds:

\begin{itemize}
\item
Case I: $S_i$'s are disjoint

\item
Case II: $|S_3\cap S_4|=i,\ i\leq k$ and $S_1,S_2,S_3\cup S_4$ are disjoint

\item
Case III: $|S_1\cap S_2|=i,\ |S_3\cap S_4|=j$ and $(S_1\cup S_2)\cap (S_3\cup S_4)=\emptyset$

\item
Case IV: $|S_2\cap S_3|=i,\ |S_3\cap S_4|=j,\ |S_2\cap S_3\cap S_4|=l,\ i,j\geq 1, l\geq 0$ and $|S_1\cap \{S_2\cup S_3\cup S_4\}|=\emptyset$

\item
Case V.: $|S_i\cap S_{i+1}|\geq 1,\ i=1,2,3$

\end{itemize}
When $S_i$'s are disjoint then the variables $X_{S_1},X_{S_2},X_{S_3},X_{S_4}$ are not independent for a fixed $n$, but if $n$ tends to infinity, then they get independent, that is $|\mathds{P}(X_{S_i}=1|X_{S_j}=1)-\mathds{P}(X_{S_i}=1)|\leq o(1)\mathds{P}(X_{S_i}=1)$.
So in Case I we can bound the sum because disjoint $S_i$'s are asymptotically independent.
In cases II-V, we can use the fact that the number of intersecting quadruples decreases as we increase the size of the intersection.

\noindent
{\bf Proof of Lemma \ref{4th_moment_bound}:}
Let us fix $T$, the tree on $k$ nodes, with $r=(r_1,\ldots r_k)$. Denote the
gluing (considering the remainder degrees $r_i$ also) of two copies of $T$ along the
set $S_i,S_j$ by $T_{i,j}=g(T_{S_i},T_{S_j})$
(if there is a valid gluing resulting in a forest).
$m_{i,j}=|V(T_{i,j})|$ and $c_{i,j}$ is just the number of components of $T_{i,j}$.
$$R_{i,j}=2(\sum_{l=1}^k r_l)-\sum_{l\in S_i\cap S_j}r_l.$$ We define similarly
$T_{i,j,h},T_{1,2,3,4},m_{i,j,h},m_{1,2,3,4},R_{i,j,h},R_{1,2,3,4}\
(i\neq j\neq h,i\neq h)$.

\begin{multline*}
\mathds{E}(A_{S_1}A_{S_2}A_{S_3}A_{S_4})=\\
=\mathds{E}(A_{S_1}A_{S_2}A_{S_3}X_{S_4})-
\mathds{E}(A_{S_1}A_{S_2}A_{S_3})\mathds{P}(X_{S_4}=1)=\\
=\mathds{P}(X_{S_4}=1)\Bigl ( \mathds{E}(A_{S_1}A_{S_2}A_{S_3}|X_{S_4}=1)-
\mathds{E}(A_{S_1}A_{S_2}A_{S_3})\Bigr )=\\
=\mathds{P}(X_{S_4}=1)\Bigl ( \mathds{P}(X_{S_3}=1)(\mathds{E}(A_{S_1}A_{S_2}|X_{S_3}=1,X_{S_4}=1)-\mathds{E}(A_{S_1}A_{S_2}|X_{S_4}=1))-\\
-\mathds{P}(X_{S_3}=1)(\mathds{E}(A_{S_1}A_{S_2}|X_{S_3}=1)-\mathds{E}(A_{S_1}A_{S_2}))
\Bigr )=\\
=\mathds{P}(X_{S_4}=1) \mathds{P}(X_{S_3}=1)\Bigl (\mathds{E}(A_{S_1}A_{S_2}|X_{S_3}=1,X_{S_4}=1)-\mathds{E}(A_{S_1}A_{S_2}|X_{S_4}=1)-\\
-\mathds{E}(A_{S_1}A_{S_2}|X_{S_3}=1)+\mathds{E}(A_{S_1}A_{S_2})
\Bigr )
\end{multline*}
Let
\begin{multline*}
f_{S_3,S_4}(X,Y)=\mathds{E}(XY|X_{S_3}=1,X_{S_4}=1)-\mathds{E}(XY|X_{S_4}=1)-\\
-\mathds{E}(XY|X_{S_3}=1)+\mathds{E}(XY),
\end{multline*}
$$f_{S_4}(X,Y,Z)=\mathds{E}(XYZ|X_{S_4}=1)-\mathds{E}(XYZ)\text{ and}$$

\begin{multline*}
\mathds{E}(A_{S_1}A_{S_2}A_{S_3}A_{S_4})=
\mathds{P}(X_{S_4}=1)\mathds{P}(X_{S_3}=1)(f_{S_3,S_4}(X_{S_1},X_{S_2})-\\
-f_{S_3,S_4}(\mathds{E}X_{S_1},X_{S_2})-
f_{S_3,S_4}(X_{S_1},\mathds{E}X_{S_2})
+f_{S_3,S_4}(\mathds{E}X_{S_1},\mathds{E}X_{S_2}))
\end{multline*}

\noindent
{\bf Case I.:} $S_i$'s are disjoint

\noindent
When $S_i$'s are disjoint, we have
$R_{1,2,3,4}=4R,\,R_{3,4}=2R,\,m_{1,2,3,4}=4k,\,m_{3,4}=2k,\,
c_{1,2,3,4}=4,\,c_{3,4}=2$.

\bigskip
Using Equation \ref{expected_value_as_polynom}
\begin{multline}
\label{ftl_2_2}
\mathds{E}(X_{S_1}X_{S_2}|X_{S_3}=1,X_{S_4}=1)
{H(r_{3,4},T_{3,4})\over H(r_{1,2,3,4},T_{1,2,3,4})}
={2^{4k}\over 3^{2R+6k-4}}\\
{n^{2R+6k-4}-n^{2R+6k-5}(6R^2+12kR-9R+18k^2-20k+{11\over 2})+O(n^{2R+6k-6})\over
n^{2R+8k-6}-n^{2R+8k-7}(2R^2+12kR-{19\over 3}R+24k^2-28k+{23\over 3})+O(n^{2R+8k-8})}
\end{multline}

\begin{multline}
\label{ftl_2_1}
\mathds{E}(X_{S_1}X_{S_2}|X_{S_3}=1)
{H(r_3,T_3)\over H(r_{1,2,3},T_{1,2,3})}=
{2^{4k}\over 3^{2R+6k-4}}\\
{n^{2R+6k-4}-n^{2R-6k-5}(4R^2+8kR-5R+12k^2-12k+{5\over 2})+O(n^{2R-6k-6})\over
n^{2R+8k-6}-n^{2R+8k-7}({4\over 3}R^2+8kR-{11\over 3}R+16k^2-16k+3)+O(n^{2R-8k-8})}
\end{multline}

\begin{multline}
\label{ftl_2_0}
{\mathds{E}(X_{S_1}X_{S_2})
\over H(r_{1,2},T_{1,2})}=
{2^{4k}\over 3^{2R+6k-4}}\\
{n^{2R+6k-4}-n^{2R+6k-5}(2R^2+4kR-1R+6k^2-4k-{1\over 2})+O(n^{2R+6k-6})\over
n^{2R+8k-6}-n^{2R+8k-7}({2\over 3}R^2+4kR-R+8k^2-4k-{5\over 3})+O(n^{2R+8k-8})}
\end{multline}

\begin{equation*}
\mathds{E}(\mathds{E}(X_{S_1})X_{S_2}|X_{S_3}=1,X_{S_4}=1)=
\mathds{P}(X_{S_1}=1)\mathds{E}(X_{S_2}|X_{S_3}=1,X_{S_4}=1)
\end{equation*}

\begin{multline}
\label{ftl_1_2}
\mathds{E}(X_{S_2}|X_{S_3}=1,X_{S_4}=1)
{H(r_{3,4},T_{3,4})\over H(r_{2,3,4},T_{2,3,4})}=
{2^{2k}\over 3^{R+3k-2}}\\
{n^{R+3k-2}-n^{R+3k-3}({5\over 2}R^2+5kR-{7\over 2}R+{25\over 2}k^2-{25\over 2}k+3)+O(n^{R+3k-4})\over
n^{R+4k-3}-n^{R+4k-4}({5\over 6}R^2+5kR-{5\over 2}R+10k^2-11k+{8\over3})+O(n^{R+4k-5})}
\end{multline}

\begin{multline}
\label{ftl_1_1}
\mathds{E}(X_{S_1}|X_{S_3}=1)
{H(r,T)\over H(r_{1,3},T_{1,3})}=
{2^{2k}\over 3^{R+3k-2}}\\
{n^{R+3k-2}-n^{R+3k-3}({3\over 2}R^2+3kR-{3\over 2}R+{15\over 2}k^2-{13\over 2}k+1)+O(n^{R+3k-4})\over
n^{R+4k-3}-n^{R+4k-4}({1\over 2}R^2+3kR-{7\over 6}R+6k^2-5k+{1\over 3})+O(n^{R+4k-5})}
\end{multline}

\begin{multline}
\label{ftl_1_0}
{\mathds{E}(X_{S_1})\over H(r,t)}=
{2^{2k}\over 3^{R+3k-2}}\\
{n^{R+3k-2}-n^{R+3k-3}({1\over 2}R^2+kR+{1\over 2}R+{5\over 2}k^2-{1\over 2}k-1)+O(n^{R+3k-4})\over
n^{R+4k-3}-n^{R+4k-4}({1\over 6}R^2+kR+{1\over 6}R+2k^2+k-2)+O(n^{R+4k-5})}
\end{multline}

\noindent
It is easy to verify that from (\ref{ftl_2_2})-(\ref{ftl_2_0}) we get that
$$f_{S_3,S_4}(X_{S_1},X_{S_2})=O(n^{-2k}),$$
and using (\ref{ftl_1_2})-(\ref{ftl_1_0})
$$
f_{S_3,S_4}(\mathds{E}(X_{S_1}),X_{S_2})=f_{S_3,S_4}(X_{S_1},\mathds{E}(X_{S_2}))=O(n^{-2k}),\ \mathrm{ and}$$
$$
f_{S_3,S_4}(\mathds{E}(X_{S_1}),\mathds{E}(X_{S_2}))=0.
$$

Hence
$\mathds{E}(A_{S_1}A_{S_2}A_{S_3}A_{S_4})=O(n^{-4k+2})$ and since the number of disjoint $S_i$ quadruples is just $n^{4k}$,
$$\sum_{S_i,\, S_i\cap S_j=\emptyset,\, i\neq j}\mathds{E}(A_{S_1}A_{S_2}A_{S_3}A_{S_4})=O(n^2).$$

\noindent
{\bf Case II.:} $|S_3\cap S_4|=i,\ i\leq k$ and $S_1, S_2, S_3\cup S_4$ are disjoint

\noindent
The number of such quadruples is $n^{4k-i}$.
Similarly as above, we want to bound $F(S_1,S_2,S_3,S_4)$.

\medskip
\noindent
From (\ref{cond_exp_diff_order}) we get
\begin{multline}
\label{2met_1}
\mathds{E}(X_{S_1}X_{S_2}|X_{S_3}=1,X_{S_4}=1)-\mathds{E}(X_{S_1}X_{S_2})=\\
=O(n^{m_{3,4}-m_{1,2,3,4}+c_{1,2,3,4}-c_{3,4}-1})=O(n^{2k-i-(4k-i)+3-1-1})=O(n^{-2k+1})
\end{multline}

\begin{multline}
\label{2met_2}
\mathds{E}(X_{S_1}X_{S_2}|X_{S_3}=1)-\mathds{E}(X_{S_1}X_{S_2})=\\
=O(n^{m_3-m_{1,2,3}+c_{1,2,3}-c_3-1})=O(n^{k-3k+3-1-1})=O(n^{-2k+1})
\end{multline}

\begin{multline}
\label{2met_3}
\mathds{E}(X_{S_1}|X_{S_3}=1,X_{S_4}=1)-\mathds{E}(X_{S_1})=\\
=O(n^{m_{3,4}-m_{1,3,4}+c_{1,3,4}-c_{3,4}-1})=O(n^{2k-i-(3k-i)+2-1-1})=O(n^{-k})
\end{multline}

\begin{multline}
\label{2met_4}
\mathds{E}(X_{S_1}|X_{S_3}=1)-\mathds{E}(X_{S_1})=\\
=O(n^{m_3-m_{1,3}+c_{1,3}-c_3-1})=O(n^{k-2k+2-1-1})=O(n^{-k})
\end{multline}

From (\ref{2met_1})-(\ref{2met_4}) it follows that
$\mathds{E}(A_{S_1}A_{S_2}A_{S_3}A_{S_4})=O(n^{-4k+3})$ and so
$$\sum_{S_i,\, |S_1\cap S_2|=i}\mathds{E}(A_{S_1}A_{S_2}A_{S_3}A_{S_4})=O(n^2).$$

\noindent
{\bf Case III.:} $|S_1\cap S_2|=i,|S_3\cap S_4|=j$ and
$(S_1\cup S_2)\cap(S_3\cup S_4)=\emptyset$.

\noindent
As before, using (\ref{cond_exp_diff_order}) we get:

\begin{multline}
\label{2_2met_1}
\mathds{E}(X_{S_1},X_{S_2}|X_{S_3}=1,X_{S_4}=1)-\mathds{E}(X_{S_1},X_{S_2})=\\
=O(n^{m_{3,4}-m_{1,2,3,4}+c_{1,2,3,4}-c_{3,4}-1})=O(n^{2k-j-(4k-i-j)+2-1-1})=O(n^{-2k+i})
\end{multline}

\begin{multline}
\label{2_2met_2}
\mathds{E}(X_{S_1},X_{S_2}|X_{S_3}=1)-\mathds{E}(X_{S_1},X_{S_2})=\\
=O(n^{m_3-m_{1,2,3}+c_{1,2,3}-c_3-1})=O(n^{k-(3k-i)+2-1-1})=O(n^{-2k+i})
\end{multline}

Now it follows from (\ref{2met_3})-(\ref{2_2met_2}) that
$\mathds{E}(A_{S_1}A_{S_2}A_{S_3}A_{S_4})=O(n^{-4k+i+2})$ and so
$$\sum_{|S_1\cap S_2|=i,|S_3\cap S_4|=j}
\mathds{E}(A_{S_1}A_{S_2}A_{S_3}A_{S_4})=O(n^2).$$

\noindent
{\bf Case IV.:} $|S_1\cap \{S_2\cup S_3\cup S_4\}|=\emptyset,|S_2\cap S_3|=i,|S_3\cap S_4|=j,
|S_2\cap S_3\cap S_4|=l,$ with $i,j\geq 1, l\geq 0$. Note, that $S_2\cup S_4\setminus S_3=\emptyset$, because otherwise there will be no valid gluing of the trees $T$ along the $S_i$'s ($g(T_{S_1}T_{S_2}T_{S_3}T_{S_4})$ is not a tree).

\noindent
The number of such $S_i$ quadruples is $n^{4k-i-j+l}$.
We want to bound
\begin{equation}
\label{bound_3met}
\mathds{E}(A_{S_1}A_{S_2}A_{S_3}A_{S_4})=
\mathds{P}(X_{S_1}=1)\left(\mathds{E}(A_{S_2}A_{S_3}A_{S_4}|X_{S_1}=1)-
\mathds{E}(A_{S_2}A_{S_3}A_{S_4})\right ).
\end{equation}
\medskip
\noindent

Using (\ref{cond_exp_diff_order}) we get:
\begin{multline}
\label{3met_3}
\mathds{E}(X_{S_2}X_{S_3}X_{S_4}|X_{S_1}=1)-
\mathds{E}(X_{S_2}X_{S_3}X_{S_4})=\\
=O(n^{m_{1}-m_{1,2,3,4}+c_{1,2,3,4}-c_{1}-1})=O(n^{k-(4k-i-j+l)+2-1-1})=O(n^{-3k+i+j-l})
\end{multline}

\begin{multline}
\label{3met_2_a}
\mathds{E}(X_{S_2}X_{S_3}|X_{S_1}=1)-
\mathds{E}(X_{S_2}X_{S_3})=\\
=O(n^{m_{1}-m_{1,2,3}+c_{1,2,3}-c_{1}-1})=O(n^{k-(3k-i)+2-1-1})=O(n^{-2k+i})
\end{multline}

\begin{multline}
\label{3met_2_b}
\mathds{E}(X_{S_3}X_{S_4}|X_{S_1}=1)-
\mathds{E}(X_{S_3}X_{S_4})=\\
=O(n^{m_{1}-m_{1,3,4}+c_{1,3,4}-c_{1}-1})=O(n^{k-(3k-j)+2-1-1})=O(n^{-2k+j})
\end{multline}

\begin{multline}
\label{3met_2_c}
\mathds{E}(X_{S_2}X_{S_4}|X_{S_1}=1)-
\mathds{E}(X_{S_2}X_{S_4})
=O(n^{m_{1}-m_{1,2,4}+c_{1,2,4}-c_{1}-1})=\\=O(n^{k-(3k-l)+c_{1,2,4}-1-1})=
\left \{
\begin{array}{ll}
O(n^{-2k+l}) & l\geq 1\\
O(n^{-2k+1}) & l=0
\end{array}
\right .
\end{multline}
Now from (\ref{2met_4}) and (\ref{3met_3})-(\ref{3met_2_c}) we get that
$$\sum_{S_1,S_2,S_3,S_4}\mathds{E}(A_{S_1}A_{S_2}A_{S_3}A_{S_4})=O(n^2),$$
where the summation is over all quadruples $S_i$ satisfying the conditions of Case IV.

\noindent
{\bf Case V.:} $|S_i\cap S_j|=s_{i,j}$ and $s_{i,i+1}\geq 1\ i=1,2,3$

\begin{equation}
\mathds{E}(X_{S_1}X_{S_2}X_{S_3}X_{S_4})=O({1\over n^{|S_1\cup S_2\cup S_3\cup S_4|-1}})
\end{equation}

\begin{equation}
\label{4_3met4}
\mathds{E}(X_{S_i}X_{S_j}X_{S_l})=O({1\over n^{|S_i\cup S_j\cup S_l|-1}})
\end{equation}

\begin{equation}
\mathds{E}(X_{S_i}X_{S_j})=O({1\over n^{|S_i\cup S_j|-c_{i,j}}})
\end{equation}
From the above it follows that
\begin{equation}
\sum_{S_1,S_2,S_3,S_4}\mathds{E}(A_{S_1}A_{S_2}A_{S_3}A_{S_4})=O(n^2),
\end{equation}
where the summation is over all quadruples $S_i$ satisfying the conditions of Case V.

\noindent
From the above the claim of Lemma \ref{4th_moment_bound} follows
\qed


\section{The limit of $X_n$}
\label{limit}

Let $U^l$ denote the set of all finite $l$-deep rooted tree.
Consider an $l$-deep rooted tree with root $x$: $T^l_x\in U^l$, with $|T^l_x|=k$. Let us denote the nodes at the $i$th level with $T_i$, and $|T_i|=t_i$ ($t_0$ is just $1$, $t_1$ is the degree of the root and $t_l$ is the number of leafs). As before $B_G(v,l)$ is the rooted $l$-ball around $v$ in $G$ and $X_n$ is a random tree on $n$ nodes with distribution $\Pi_n$.

If we assign remainder degrees $r$ to the rooted tree $T_x^l$ and forget the root, then using Lemma \ref{forest_prob} we get
\begin{equation}
\label{exp_subtree}
\mathds{E}(X_n^T)=\mathds{E}\left(\sum_{S\subseteq [n]}X_S^T\right)={n!\over (n-k)!}\mathds{P}(X_S^T=1)=n{n-1\over n-k}{{3n-R-3k-1\choose n-R-k}\over {3n-3\choose n-2}}H(r,T).
\end{equation}
Now the number of vertices $v\in X_n$ for which $B_{X_n}(v,l)\cong T_x^l$ is just
\begin{equation}
\label{neigh_stat}
{1\over |Aut(T^l_x)|}
\sum_{r_k=1}^{n}\sum_{r_{k-1}=1}^{n}\cdots \sum_{r_{k-t_l+1}=1}^{n} X_n^T
\end{equation}
Where the remainder degrees are $0$ except for the leafs of $T^l_x$ ($r_i=0,\ \forall i\notin T_l$)
and $Aut(T^l_x)$ means the set of all rooted automorphisms of $T^l_x$.

\begin{remark}
\label{max_deg_rem}
Using lemma \ref{lemma_maximum_also_becsles} we have that the maximum degree is almost surely less than $3\log_3n$. For a tree $T$ with remainder degree $r$, if there exists
$i$, s.t. $r_i>3\log_3n$, then $X_n^T=0$ almost surely. So we can omit the terms with $r_i>3\log_3n$ in \ref{neigh_stat}.
\end{remark}
Now if $v$ is a uniform random vertex of $X_n$, then using remark \ref{max_deg_rem} we have

\begin{equation}
\label{labsubg_neigh}
\mathds{P}(B_{X_n}(v,l)=T_x^l|X_n)={1\over |Aut(T^l_x)|}
\sum_{r_k=1}^{3\log_3n}\sum_{r_{k-1}=1}^{3\log_3n}\cdots \sum_{r_{k-t_l+1}=1}^{3\log_3n} {X_n^T\over n},
\end{equation}
Now from Theorem \ref{labeled_tree_deviation}, we get that
\begin{multline*}
\mathds{P}(B_{X_n}(v,l)=T_x^l)={1\over |Aut(T^l_x)|}
\sum_{r_k=1}^{3\log_3n}\sum_{r_{k-1}=1}^{3\log_3n}\cdots \sum_{r_{k-t_l+1}=1}^{3\log_3n} \left({\mathds{E}(X_n^T)\over n}+o({1\over n^{{1\over 8}}})\right)=\\
=o(1)+{1\over |Aut(T^l_x)|}
\sum_{r_k=1}^{3\log_3n}\sum_{r_{k-1}=1}^{3\log_3n}\cdots \sum_{r_{k-t_l+1}=1}^{3\log_3n} {\mathds{E}(X_n^T)\over n}.
\end{multline*}

$$
\lim_{n\rightarrow \infty}{n-1\over n-k}
{{3n-R-3k-1\choose n-R-k}\over {3n-3\choose n-2}}
H(r,T)=
9\left({4\over 27}\right)^{\sum_0^l t_j}{\sum_{i\in T_l} r_i \prod_{i\in T_l} (r_i+1) \over 3^{\sum_{i\in T_l} r_i}}\prod_{j\notin T_l} d_j!,
$$
so it follows that
\begin{multline}
\label{negdens}
\lim_{n\rightarrow \infty}\mathds{P}(B_{X_n}(v,l)=T_x^l)=\\
=9\left({4\over 27}\right)^{\sum_0^{l-1} t_j}{t_l\over 3^{t_l}}{1\over
|Aut(T_x^l)|}\prod_{j\notin T_l}d_j!=p(T_x^l).
\end{multline}

Let $\cal{G}$ denote the set of all countable connected rooted trees. Let ${\cal T}({\cal G},T_x^l):=\{(G,v)\in {\cal G}:B_G(v,l)\cong T_x^l\}$. Now consider a measure $\mu$ on the sets ${\cal T}({\cal G},T_x^l)$:
$$
\mu ({\cal T}({\cal G},T_x^l)=p(T_x^l)
$$

\begin{lemma}
\label{mu_prob_measure}
$\mu$ extends to a probability measure on $\cal G$.
\end{lemma}

\noindent
{\bf Proof:} We only need to show that
\begin{equation}
\label{prob_meas}
p(T_x^{l-1})=\sum_{\displaystyle T_x^{l}: B_{T_x^{l}}(x,l)\cong T_x^{l-1}}p(T_x^{l})
\end{equation}

\noindent
Let $\sigma, \rho\in Aut(T_x^l)$ be two automorphisms of the rooted tree $T_x^l$. We say that $\sigma \sim \rho$ if and only if there exists $\tau \in Aut(T_x^l)$, such that $\tau$ fixes every vertex not in $T_{l}$ and $\sigma \circ \tau = \rho$. $\sim$ is an equivalence. The equivalence classes have $\prod_{i\in T_{l}}(d_i-1)!$ elements, hence it follows
$$|Aut(T_x^l)|=|Aut(T_x^l\setminus T_{l})|\prod_{i\in T_{l}}(d_i-1)!.$$
Now we get equation \ref{prob_meas} easily from the following:
\begin{multline*}
\label{lab_dist_sum}
\sum_{d_{n-t_l-t_{l-1}+1}=1}^{\infty}\sum_{d_{n-t_l-t_{l-1}+2}=1}^{\infty}\cdots
\sum_{d_{n-t_l}=1}^{\infty}
{9\left({4\over 27}\right)^{\sum_0^{l-1} t_j} {t_l\over 3^{t_l}}\prod_{j\notin T_l}d_j!\over |Aut(T_x^l)|}=\\
={9\left({4\over 27}\right)^{\sum_0^{l-1} t_j} \left({9\over 4}\right)^{t_{l-1}}t_{l-1}\prod_{j\notin T_l\cup T_{l-1}}d_j!\over
|Aut(T_x^l\setminus T_l)|}.
\end{multline*}
\qed

From Lemma \ref{mu_prob_measure} it follows that $\mu$ is indeed the limit of the random tree sequence $X_n$. This completes the proof of Theorem \ref{tree_limit}.

\end{document}